\documentclass[12pt,twoside]{article}
\usepackage[T1]{fontenc}
\usepackage{textcomp}
\usepackage{mathrsfs}
\usepackage{amsmath}
\usepackage{amssymb,bm}
\usepackage{fancyhdr}
\usepackage{latexsym}
\usepackage{bbding}
\usepackage{mathrsfs}
\usepackage{wasysym}
\usepackage{cite}
\usepackage{hyperref}
\usepackage{multicol,graphics}
\usepackage{xcolor}

\RequirePackage{times}

\def\XXint#1#2#3{{\setbox0=\hbox{$#1{#2#3}{\int}$ }
		\vcenter{\hbox{$#2#3$ }}\kern-.6\wd0}}

\newtheorem{theorem}{Theorem}[section]
\newtheorem{lemma}[theorem]{Lemma}
\newtheorem{corollary}[theorem]{Corollary}

\newtheorem{definition}[theorem]{Definition}
\newtheorem{remark}[theorem]{Remark}

\numberwithin{equation}{section}
\newenvironment{proof}[1][Proof]{\noindent\textbf{#1.} }{\hfill $\Box$}

\allowdisplaybreaks[4]
\makeatletter
\begin{document}
	
	\title{Energy equality of the weak solutions to the fractional  Navier-Stokes  / MHD equations\footnote{Authors are listed alphabetically by surname then given name. Authors equally share the first authorship. }}
	\author{ Yi Feng\footnote{mx120220310@stu.yzu.edu.cn,fengyimhd@163.com}, Weihua Wang\footnote{Corresponding author: wangvh@163.com, wangweihua15@mails.ucas.ac.cn (W.Wang).}\\
		[0.2cm] {\small School of Mathematical Science, Yangzhou University,}\\
		[0.2cm] {\small  Yangzhou, Jiangsu 225009,  China}}
\date{}	
	\maketitle
	\begin{abstract}
		In this paper, we study the problem of energy equality for weak solutions of the 3D incompressible fractional   Navier-Stokes / MHD equations. With the help of  the technique of symmetrization and interpolation method, we obtain some new sufficient conditions including the Sobolev multiplier spaces, which insures the validity of the energy equality of the weak solution to fractional MHD equations. Correspondingly, the results of fractional Navier-Stokes equations are obtained.   And these energy equations are usually related to the uniqueness of solutions to the corresponding fractional   Navier-Stokes / MHD equations.
	\end{abstract}
	\smallbreak
	
	\textit{Keywords}: fractional MHD,  fractional  Navier-Stokes equations,  Energy equality,  Onsager's conjecture,  Multiplier spaces.
	\smallskip
	
	\textit{2020 AMS Subject Classification}:  76W05; 76B03; 35Q35; 76D05
	
\section{Introduction}\label{sec1}

We consider the issue of energy equality for the weak solutions to the following fractional  Navier-Stokes  / MHD equations with initial  value in the appropriate space
	\begin{equation}\label{eq1.1}
		\begin{cases}
			u_{t}+(-\triangle)^{\alpha} u +(u\cdot \nabla)u + \nabla p=0,\quad &(x,t)\in \mathbb{R}^{3}\times (0,T),\\
			\text{div} u =0, \quad &(x,t)\in \mathbb{R}^{3}\times (0,T),\\
			u(x,0)=u_{0},  \quad &x\in \mathbb{R}^{3}
		\end{cases}
	\end{equation}
	and
	\begin{equation}\label{eq1.2}
		\begin{cases}
			u_{t}+(-\triangle)^{\alpha} u +(u\cdot \nabla)u - (B\cdot \nabla)B + \nabla p=0,\quad &(x,t)\in \mathbb{R}^{3}\times (0,T),\\
			B_{t}+(-\triangle)^{\alpha} B +(u\cdot \nabla)B - (B\cdot \nabla)u =0,\quad &(x,t)\in \mathbb{R}^{3}\times (0,T),\\
			\text{div} u = \text{div} B =0, \quad &(x,t)\in \mathbb{R}^{3}\times (0,T),\\
			u(x,0)=u_{0}, \quad B(x,0)=B_{0}, \quad &x\in \mathbb{R}^{3},
		\end{cases}
	\end{equation}
	where $u(x,t),B(x,t)$ and $p(x,t)$ represent the velocity fields, magnetic fields, and scalar pressure, respectively. $u_0$ and $B_0$ denote initial value of velocity fields and magnetic fields, which also satisfying $\text{div}u_0=\text{div}B_0=0$ in distribution sense.  The fractional operator $(-\triangle)^{\alpha}$ is defined at the Fourier level by the symbol $|\xi|^{2{\alpha}}$. And we also write  $\Lambda^{\alpha} \triangleq (-\triangle)^{\frac{\alpha}{2}}$ for any $\alpha\in \mathbb{R}$.
	
	In the case of $\alpha=1$,  Leray \cite{Leray34} and Hopf \cite{Hopf51} showed the global existence of a weak solution to \eqref{eq1.1}($\alpha=1$) with the energy inequality (opposed to the equality \eqref{eq1.4})
 \begin{equation*}
  \frac{1}{2}\|u(t)\|^{2}_{L^{2}} + \int_{0}^{t}\|\nabla u(\tau)\|^{2}_{L^{2}}{\rm d}\tau\leq \frac{1}{2}\|u_0\|^{2}_{L^{2}},\quad t \geq 0.
  \end{equation*}
  The above energy inequality cannot be equated because the weak solution lacks regularity and it comes from a limiting process on the smoother functions based on weak convergence.
	 The problem of the energy equality is an important property of the weak solutions of Navier-Stokes equation, MHD equation and many other fluid mechanics equations.  And these energy equalities are linked to the Onsager's conjecture \cite{Onsager1949} of the corresponding  fluid mechanics equations.

In particular, for the Navier-Stokes equations, Lions \cite{LJL60} obtained energy equality  for the Leray-Hopf weak solutions in $L^{4}(0,T; L^{4})$.  Whereafter, Serrin \cite{SJ63} proved that  the energy equality of the Leray-Hopf weak solutions holds in $L^{q}(0,T;L^{p}(\Omega))$ with $\frac{2}{q} + \frac{d}{p} \leq 1$. And Shinbrot \cite{SM74} optimized the condition to $\frac{2}{p}+\frac{2}{q}\leq 1$ with $p\ge4$, which does not depend on the space dimension $d$.  Galdi \cite{GGP19} established the energy equality of the distributional solutions in $\in L^{4}(0,T;L^{4}(\Omega))$  by mollifying procedure and duality argument, which improved Lion's result \cite{LJL60}.  Berselli and Chiodaroli \cite{LCE20} obtained some new criteria for the weak solutions, involving $\nabla u$, and proved that if $u\in L^{q}(0,T;L^{p}(\Omega))$ with $\frac{2}{q}+\frac{2}{p}= 1$, the energy equality of the distributional solutions to Navier-Stokes equations is true, which follows by a duality argument as done in \cite{GGP19} and involves the full range of exponents. Beir\~{a}o and Yang \cite{HJ19} concerned with the energy equality for the weak solutions to the non-Newtonian equations, and presented some conditions that contain $\nabla u$, which included the conclusions in Berselli and Chiodaroli \cite{LCE20}. Cheskidov and  Luo \cite{ChL20} established the energy equality for weak solutions in a large class of function spaces including the critical Onsager space $L^{3}B^{\frac{1}{3}}_{3,\infty}$. And Chen and Zhang \cite{CZ23} gave some sufficient conditions, including the multiplier spaces,  of the energy equality holding. Very recently, Wu \cite{Wu24} obtained the energy equality for the distributional solutions to \eqref{eq1.1} in addition to the  conditions $u$ in $L^{\frac{4\alpha}{2\alpha-1}}(0,T;L^{4}(\mathbb{R}^{3}))$ with $\alpha\geq 1$.

	For the MHD equations, we  need to possess more complicated nonlinear terms. He and Xin \cite{HeX05} established the validity of the energy equality for $u$ in the Serrin class (i.e. $L^{q}(0,T;L^{p}(\Omega))$ with $\frac{2}{q} + \frac{3}{p} = 1, ~p\geq3$) even though $B$  in
$L^{\infty}(0,T; L^{2}(\mathbb{R}^{3}))\cap L^{2}(0,T; H^{1}(\mathbb{R}^{3}))$ given by Leray \cite{Leray34}  and Hopf \cite{Hopf51}.
 Under the condition that $u$ belongs to $L^{\infty}(0,T;L^{n}(\Omega))$,  Yan and Jiu \cite{YJ09} obtained the energy equality of the weak solutions without the additional assumption of magnetic field.   Zeng \cite{ZY22} proved the energy equality of the weak solution in $L^{q}(0,T;L^{p}(\Omega))$ with : (\romannumeral1) $\frac{1}{p}+\frac{1}{q}\leq \frac{1}{2}, ~p\geq4$ or (\romannumeral2) $\frac{3}{p}+\frac{1}{q}\leq 1, p\in (3,4]$, which cancels the condition $\frac{2}{p}+\frac{2}{q}\leq 1$ of Taniuchi\cite{TY97}.
By Galerkin technique, Lai and Yang \cite{BY22} obtained the validity of   the energy equality for the distributional solutions in $ L^{4}(0,T;L^{4}(\mathbb{R}^{3}))$.  Tan and Wu \cite{TW22} established the energy equality of weak solution in $L^{2, \infty}(0,T; BMO_{loc}(\mathbb{R}^{3}))$.   Wang and Zuo \cite{WangZ22} dealt with the critical case of Shinbrot \cite{SM74} and the corresponding result to He and Xin \cite{HeX05}  in the bounded
domain.   Eguchi \cite{Eguchi23} established the energy equality   for a larger class of the magnetic field  $B$ than that of the velocity field $u$ in a bounded domain, which is as a generalization of He and Xin \cite{HeX05}.

The fractional Navier-Stokes equations \eqref{eq1.1} came from  Lions \cite{LJL69} has a unique global smooth solution for any smooth initial data  provided $\alpha\geq\frac{5}{4}$ (also see \cite{WuJH03}).  However, when $\alpha<\frac{5}{4}$, whether the fractional Navier-Stokes equation \eqref{eq1.1} with large initial values has a smooth solution remains a well-known open problem, just like the classical Navier-Stokes equation.  Similar to the research strategy of energy conservation in classical Navier-Stokes equations, {\bf our aim} is to establish the energy equality for a larger class of solutions for the fractional Navier-Stokes / MHD equations. In this paper, inspired by Chen and Zhang \cite{CZ23}, with the help of symmetrization, we obtain some sufficient conditions, containing the multiplier spaces,  which guarantee the validity of the energy equality of the weak solutions to fractional  MHD equations.  Moreover, we also give the corresponding results of fractional Navier-Stokes equation.
 We extend  Chen and Zhang's  results \cite{CZ23} to the  fractional Navier-Stokes / MHD equations. And we also generalize Wu's results in \cite{Wu24} and Zeng's results in \cite{ZY22}.

 To accurately describe our result, we first recall  the definition of the weak solution to \eqref{eq1.1} and \eqref{eq1.2}.
 \begin{definition}[Leray-Hopf Weak Solutions]\label{de2.5} We say that
		 $u\in L^{\infty}(0,T;L^{2}_{\sigma}(\mathbb{R}^{3}))\cap L^{2}(0,T;H^{\alpha}_{0}(\mathbb{R}^{3}))$ is a Leray-Hopf weak solution to \eqref{eq1.1} with initial values $u_0$ in $L^{2}_{\sigma}(\mathbb{R}^{3})$, if  the following conditions holds:
\begin{description}
		\item[(\romannumeral1)]  $u$ is a solution of \eqref{eq1.1} in the sense of distribution
	                             \begin{equation*}
			                      \int_{0}^{T} \left((\psi_\tau + (u\cdot\nabla)\psi, u) -(u, \Lambda^{2\alpha}\psi) \right){\rm d}\tau = -(u_0,\phi(0)),
		                        \end{equation*}
		                        for all $\psi \in C^{\infty}_0([0,T);C^{\infty}_{0,\sigma}(\mathbb{R}^{3}))$;
		\item[(\romannumeral2)]  $u$ satisfies the energy inequality
                                 \begin{equation*}
        	                        \frac{1}{2}\|u(t)\|^{2}_{L^{2}} + \int_{0}^{t}\|\Lambda^{\alpha} u(\tau)\|^{2}_{L^{2}}{\rm d}\tau\leq  \frac{1}{2}\|u_0\|^{2}_{L^{2}},\quad t \in (0,T);
                                 \end{equation*}
		\item[(\romannumeral3)]  $\|u(t)-u_0\|_{L^{2}}\to 0 \text{~~as~~} t\to 0^{+}$.
	\end{description}	
	\end{definition}

		\begin{definition}[Leray-Hopf Weak Solutions]\label{de2.6}
		$\begin{pmatrix}
			u\\
			B
		\end{pmatrix} \text{ in } L^{\infty}(0,T;L^{2}_\sigma(\mathbb{R}^{3}))\cap L^{2}(0,T;H^{\alpha}_{0}(\mathbb{R}^{3}))$  is called  Leray-Hopf weak solution of equation \eqref{eq1.2}  with initial values $u_0$ and  $B_0$ in $L^{2}_{\sigma}(\mathbb{R}^{3})$, if $\begin{pmatrix}
			u\\
			B
		\end{pmatrix}$ satisfies
	\begin{description}
		\item[(\romannumeral1)]  $\begin{pmatrix}
			u\\
			B
		\end{pmatrix}$ is a solution of \eqref{eq1.2} in the sense of distribution
\begin{eqnarray*}
 \int_{0}^{T} \left((\phi_\tau + (u\cdot\nabla)\phi, u) -(u, \Lambda^{2\alpha}\phi) -((B\cdot\nabla)\phi, B)\right){\rm d}\tau &=& -(u_0,\phi(0)), \\
  \int_{0}^{T} \left((\psi_\tau + (u\cdot\nabla)\psi, B) -(B, \Lambda^{2\alpha}\psi) -((B\cdot\nabla)\psi, u)\right){\rm d}\tau &=& -(B_0,\psi(0))\\
\end{eqnarray*}
for all $\phi, \psi \in C^{\infty}_0([0,T);C^{\infty}_{0,\sigma}(\mathbb{R}^{3}))$;
		\item[(\romannumeral2)]  $\begin{pmatrix}
			u\\
			B
		\end{pmatrix} $  satisfies the energy inequality
		\begin{equation*}
			 \frac{1}{2}\|u(t)\|^{2}_{L^{2}} +  \frac{1}{2}\|B(t)\|^{2}_{L^{2}} + \int_{0}^{t}\left(\|\Lambda^{\alpha} u(\tau)\|^{2}_{L^{2}}+\|\Lambda^{\alpha} B(\tau)\|^{2}_{L^{2}}\right){\rm d}\tau\leq \frac{1}{2}\|u_0\|^{2}_{L^{2}} +  \frac{1}{2}\|B_0\|^{2}_{L^{2}}
		\end{equation*}
		for any $t \in (0,T)$.
		\item[(\romannumeral3)]  $\left\|\begin{pmatrix}
			u(t)\\
			B(t)
		\end{pmatrix} - \begin{pmatrix}
			u_{0}\\
			B_{0}
		\end{pmatrix}\right\|_{L^{2}}\to 0 \text{ ~as~ } t\to 0^{+}$.	
	\end{description}		
\end{definition}

 Specifically, our results are as follows:
	\begin{theorem}\label{th1.1}
	Let $\begin{pmatrix}
		u\\
		B
	\end{pmatrix}$  be a  weak solutions of the fractional MHD equations \eqref{eq1.2}   on some time interval $[0,T]$ with $u_0,B_0\in L^{2}_{\sigma}(\mathbb{R}^{3})$ and $0<T\leq\infty$. If $u$ and $B$ satisfy $|u|^{\frac{2}{\beta+\theta}-2}u$ and $|B|^{\frac{2}{\beta+\theta}-2}B \in L^{\frac{2\alpha(\beta+\theta)}{\alpha-\theta}}\mathcal{M}^{3}(\dot{H}^{\frac{\theta}{\beta+\theta}}\to L^{2(\beta+\theta)})$, where $\alpha\geq1,~\beta,\theta\ge0,\,  \frac{1}{2}\le\beta+\theta\le1$,$\,$
	then the following energy equality
	\begin{equation}\label{eq1.3}
		\frac{1}{2}\|u(t)\|^{2}_{L^{2}} + \frac{1}{2}\|B(t)\|^{2}_{L^{2}} + \int_{0}^{t}\left(\|\Lambda^{\alpha} u(\tau)\|^{2}_{L^{2}}+\|\Lambda^{\alpha} B(\tau)\|^{2}_{L^{2}}\right){\rm d}\tau = \frac{1}{2}\|u_0\|^{2}_{L^{2}} + \frac{1}{2}\|B_0\|^{2}_{L^{2}} .
	\end{equation}
holds for any $\;0\le t<T$.
\end{theorem}
\begin{remark}
  Comparing the  results in \cite{ZY22}, taking $\alpha=1$ in Theorem \ref{th1.1}, obviously, we can observe that $\frac{2(\beta+\theta)}{1-\theta}\in[1,+\infty)$, which is greater than the condition $q\geq 4$ in \cite{ZY22}. Considering the following embedding relationship,	$$L^{\frac{6(\beta+\theta)}{3-3\beta-\theta}}(\mathbb{R}^{3})\hookrightarrow\mathcal{M}^{3}(\dot{H}^{\frac{\theta}{\beta+\theta}}(\mathbb{R}^{3})\to L^{2(\beta+\theta)}(\mathbb{R}^{3})),$$
   $\frac{6(\beta+\theta)}{3-3\beta-\theta}\in [\frac{6}{5},+\infty)$ in Theorem \ref{th1.1}  optimizes the condition $p\in (3,4]$ in \cite{ZY22}.
\end{remark}
	
\begin{corollary}\label{co1.2}
		Let $u$ be a weak solution of the fractional Navier-Stokes equations \eqref{eq1.1}   on some time interval $[0,T]$ with $u_0\in L_{\sigma}^{2}(\mathbb{R}^{3})$ and $0<T\leq\infty$. If $u$ satisfies $|u|^{\frac{2}{\beta+\theta}-2}u\in L^{\frac{2\alpha(\beta+\theta)}{\alpha-\theta}}\mathcal{M}^{3}(\dot{H}^{\frac{\theta}{\beta+\theta}}\to L^{2(\beta+\theta)})$, where $ \alpha\geq1,~\beta,~\theta\ge0,~  \frac{1}{2}\le\beta+\theta\le1$, the energy equality
\begin{equation}\label{eq1.4}
		\frac{1}{2}\|u(t)\|^{2}_{L^{2}} +  \int_{0}^{t}\|\Lambda^{\alpha} u(\tau)\|^{2}_{L^{2}}{\rm d}\tau = \frac{1}{2}\|u_0\|^{2}_{L^{2}}.
	\end{equation}
 holds for any $0\le t<T$.
\end{corollary}

	\begin{theorem}\label{th1.3}
		Let $\begin{pmatrix}
			u\\
			B
		\end{pmatrix}$  be a  weak solutions of the fractional MHD equations \eqref{eq1.2}   on some time interval $[0,T]$ with $u_0,B_0\in L^{2}_{\sigma}(\mathbb{R}^{3})$ and $0<T\leq\infty$. If $u$ and $B$ satisfy
		\begin{equation*}
			|u|^{\frac{10\theta^{2}-23\theta+7}{3+8\theta-5\theta^{2}}}u,~|B|^{\frac{10\theta^{2}-23\theta+7}{3+8\theta-5\theta^{2}}}B \in L^{\frac{2\alpha(-5\theta^{2}+8\theta+3)}{2\alpha(-\theta+4)+5\theta^{2}-8\theta-3}}\mathcal{M}^{3}(\dot{H}^{1}\to L^{\frac{-5\theta^{2}+8\theta+3}{3-2\theta}})
		\end{equation*}
and
		\begin{equation*}
			|\nabla u|^{\frac{10-5\theta}{7-3\theta}},~|\nabla B|^{\frac{10-5\theta}{7-3\theta}} \in L^{\frac{2\alpha(7-3\theta)(\beta+1)}{4\alpha(10-5\theta)(\beta+1)-(7-3\theta)(4\alpha\beta+4\alpha+1)}}\mathcal{M}^{3}(\dot{H}^{\frac{1}{\beta+1}}\to L^{1}),
		\end{equation*}
		 where $ \alpha\geq1,~\beta\ge0,~0\leq\theta<\frac{1}{2}$, the energy equality \eqref{eq1.3} holds for any $0\le t<T$.
		
	\end{theorem}
	\begin{corollary}\label{co1.4}
			Let $u$ be a weak solution of the fractional Navier-Stokes equations \eqref{eq1.1}   on some time interval $[0,T]$ with $u_0\in L_{\sigma}^{2}(\mathbb{R}^{3})$ and $0<T\leq\infty$. If $u$ satisfies
		\begin{equation*}
		|u|^{\frac{10\theta^{2}-23\theta+7}{3+8\theta-5\theta^{2}}}u\in L^{\frac{2\alpha(-5\theta^{2}+8\theta+3)}{2\alpha(-\theta+4)+5\theta^{2}-8\theta-3}}\mathcal{M}^{3}(\dot{H}^{1}\to L^{\frac{-5\theta^{2}+8\theta+3}{3-2\theta}})
		\end{equation*}
and
		\begin{equation*}
		|\nabla u|^{\frac{10-5\theta}{7-3\theta}} \in L^{\frac{2\alpha(7-3\theta)(\beta+1)}{4\alpha(10-5\theta)(\beta+1)-(7-3\theta)(4\alpha\beta+4\alpha+1)}}\mathcal{M}^{3}(\dot{H}^{\frac{1}{\beta+1}}\to L^{1}),
		\end{equation*}
		where $ \alpha\geq1,~\beta\ge0,~0\leq\theta<\frac{1}{2}$, the energy equality \eqref{eq1.4} holds for any $0\le t<T$.
\begin{remark}
 The assumptions of Corollary \ref{co1.2} and Corollary \ref{co1.4} cover two assumptions in \cite{CZ23}.
\end{remark}		
	\end{corollary}
	\begin{theorem}\label{th1.5}
		Let $\begin{pmatrix}
			u\\
			B
		\end{pmatrix} $  be a weak solutions of the fractional MHD equations \eqref{eq1.2}   on some time interval $[0,T]$ with $u_0,B_0\in L^{2}_{\sigma}(\mathbb{R}^{3})$ and $0<T\leq\infty$. If $u$ and $B$ satisfy
		\begin{equation*}
			\nabla u,~\nabla B \in L^{p}\mathcal{M}^{3}(\dot{W}^{\frac{1+\theta}{1+\beta},\frac{6-3\theta}{3-\theta^{2}}}\to L^{\frac{12-6\theta}{9-5\theta}})
		\end{equation*}
and
		\begin{equation*}
			 u, ~B \in L^{q}(L^{\frac{12-6\theta}{3-\theta}}),~ \frac{1}{p}+\frac{1}{q}=\frac{1}{2},~2\leq p<\infty,
		\end{equation*}
		where $ \alpha\geq1,~\beta\ge0,~-1\leq\theta<0$, the energy equality \eqref{eq1.3} holds for any $0\le t<T$.
		
	\end{theorem}
	\begin{corollary}\label{co1.6}
		Let $u$ be a weak solution of the fractional Navier-Stokes equations \eqref{eq1.1}   on some time interval $[0,T]$ with $u_0\in L_{\sigma}^{2}(\mathbb{R}^{3})$ and $0<T\leq\infty$. If $u$ satisfies
		\begin{equation*}
			\nabla u \in L^{p}\mathcal{M}^{3}(\dot{W}^{\frac{1+\theta}{1+\beta},\frac{6-3\theta}{3-\theta^{2}}}\to L^{\frac{12-6\theta}{9-5\theta}})
		\end{equation*}
and
		\begin{equation*}
			u \in L^{q}(L^{\frac{12-6\theta}{3-\theta}}),~ \frac{1}{p}+\frac{1}{q}=\frac{1}{2},~2\leq p<\infty,
		\end{equation*}
		where $ \alpha\geq1,~\beta\ge0,~-1\leq\theta<0$, the energy equality \eqref{eq1.4} holds for any $0\le t<T$.
		
	\end{corollary}
\begin{remark}
   The range of temporal integrability exponent $q$ in Corollary \ref{co1.6} is $2\leq q<\infty$, which extends the range of temporal integrability exponent $\frac{4-2\theta}{1-\theta}\in [3,4)$ in (\romannumeral3) of {\rm\cite[Theorem 1]{CZ23}}.
\end{remark}

	The remainder of this paper is organized as follows: In Section \ref{Sec2}, we outline some basic definitions and facts, as well as give some key lemmas.     And in Section \ref{Sec3}, we give the proofs of Theorem \ref{th1.1} and Corollary \ref{co1.2}. In Section \ref{Sec4}, we prove the Theorem \ref{th1.3} and Theorem \ref{th1.5}, Corollary \ref{co1.4} and Corollary \ref{co1.6}.\\

\noindent {\bf Notation}.  ~~Throughout the paper,  $A \lesssim B$ denotes $|A|\le C|B|$ with some positive constant $C$. If $C$ depends on a parameter $\varepsilon$,  that is, $C=C(\varepsilon)$, ~$\lesssim$ is replaced by $ \lesssim_{\varepsilon}$.

\section{Preliminaries}\label{Sec2}
In this section, 	we will recall some function spaces,  mollifier, multiplier spaces and Gagliardo-Nirenberg inequality. And we also recall the definitions of weak solution for the fractional Navier-Stokes  / MHD equations and some remarks.

Let $\Omega$ is  a region of $\mathbb{R}^{3}$ or $\mathbb{R}^{3}$.  The Lebesgue space $L^{p}(\Omega)$ consists of all strongly measurable functions $u: \Omega\to \mathbb{R}^d$ with the norm
	\begin{equation*}
		\|u\|_{L^{p}(\Omega)}:=\begin{cases}
  \left(\int_{\mathbb{R}^{3}}|u(x)|^{p}dx\right)^{\frac{1}{p}}, &1\leq p<\infty;\\
  \underset{\Omega}{\text{ess}\sup}|u(x)|, & p=+\infty.
\end{cases}
	\end{equation*}
The space $L^{q}(0,T;~L^{p}(\Omega))$ consists of all strongly measurable functions $u: [0,T]\to  L^{p}(\Omega) $ with	
\begin{equation*}
		\|u\|_{L^{q}(0,T;L^{p}(\Omega)}
:=\begin{cases}
  \left(\int_{0}^{T}\|u(t)\|_{L^{p}(\Omega)}^{q}dt\right)^{\frac{1}{q}}, & 1\leq q<\infty;\\
   \underset{0\leq t\leq T}{\text{ess}\sup}\|u(t)\|_{L^{p}(\Omega)}, & q=+\infty.
  \end{cases}
\end{equation*}

	For convenience, we denotes $\|\cdot\|_{L^{q}(L^{p})}$ by the norm $\|\cdot\|_{L^{q}(0,T;L^{p}(\mathbb{R}^{3})}.\;$

	The function spaces of vector fields used in this paper are as follows:
	\begin{equation*}
		\begin{split}
			&C^{\infty}_{c,\sigma} (\mathbb{R}^{3})=\{u\in C^{\infty}_c (\mathbb{R}^{3}),\text{div} ~u=0\},\\
			&C^{\infty}_{0,\sigma} (\mathbb{R}^{3})=\{u\in C^{\infty}_0 (\mathbb{R}^{3}),\text{div} ~u=0\},\\
			&H^{\alpha}_{0,\sigma} (\mathbb{R}^{3})=\{u\in H^{\alpha} (\mathbb{R}^{3}),\text{div}~u=0\},
		\end{split}
	\end{equation*}

	Let $\eta$ be a standard mollifier \cite[\S.5, ~p713]{ELJ98}, that is, $\eta(t)=C\exp\left(\frac{1}{|t|^{2}-1}\right)\chi_{\{|t|\leq 1\}}$, where the constant $C\geq 0$ selected to integrate to unity and $\chi_{\{|t|\leq 1\}}$ is the indicator function of $\{|t|\leq 1\}$.  For any $\varepsilon>0$,  we set the rescaled mollifier $\eta_{\varepsilon}(t):=\frac{1}{\varepsilon}\eta(\frac{t}{\varepsilon})$, then $\int_{\mathbb{R}}\eta_{\varepsilon}(t){\rm d}t=1$.  Let $s, ~t$ be two fixed  numbers such that $0\leq s < t<\infty$.  For any function $f\in L^{1}_{loc}$, define its mollification
	\begin{equation*}
		f_{\varepsilon}(\tau)=J_{\varepsilon}[f](\tau)=\int_{s}^{t}\eta_\varepsilon(\tau-\sigma)f(\sigma){\rm d}\sigma.
	\end{equation*}
	
For the reader's convenience, we give the following lemma for the mollifier.
   \begin{lemma}\cite[Lemma 2.1]{Masuda84}\label{Lem2.0}
    Let $X$ be a Banach space and any $u$ be in $L^{q}((s,t); X), ~1\leq q <+\infty$. We have
    \begin{description}
      \item[(\romannumeral1)] For each fixed $\varepsilon$, $J_{\varepsilon}$ is a bounded operator from $L^{q}((s,t); X)$ into $C^{1}([s,t]; X)$;
      \item[(\romannumeral2)]  For each fixed $u$ in $L^{q}((s,t); X)$, $J_{\varepsilon}[u]\to u$ as $\varepsilon\to 0$ in $L^{q}((s,t); X)$;
      \item[(\romannumeral3)]  If $u$ be in $C^{1}([s,t]; X)$, then $J_{\varepsilon}[u](t)\to \frac{1}{2}u(t)$ and $J_{\varepsilon}[u](s)\to \frac{1}{2}u(s)$ as $\varepsilon\to 0$ in the norm of $X$.
    \end{description}
   \end{lemma}

	Next we brief  the definitions on multiplier spaces  and related comments.
	\begin{definition}[\cite{GO09}]\label{de2.1}
		Let $E$ and $F$ be two Banach spaces embedded in $S'(\mathbb{R}^{3})$.  The multiplier spaces
\begin{equation*}
  \mathcal{M}(E\to F) := \left\{ \begin{aligned}\text{the collection of all functions~} \varphi \text{~satisfying~}  \varphi u\in F \\ \text{~for any~}u\in E  \text{~and the norm~}  \|\varphi\|_{\mathcal{M}(E\to F)} <\infty\end{aligned} \right\}
\end{equation*}
where
        \begin{equation*}
			\|\varphi\|_{\mathcal{M}(E\to F)} \triangleq \sup_{u\in E\setminus\{0\}}\frac{\|\varphi u\|_F}{\|u\|_E}.
		\end{equation*}
		Furthermore, if $\bm{\varphi}$ is a $d$-dimension vector-value function, i.e. $\bm{\varphi} = (\varphi^{1},...,\varphi^{d})\in\;\mathcal{M}(E\to F),\;$the norm becomes
		\begin{equation*}
			\|\bm{\varphi}\|^{2}_{\mathcal{M}^{d}(E\to F)}:=\sum_{1\le i\le d}	\|\varphi^{i}\|^{2}_{\mathcal{M}(E\to F)}.
		\end{equation*}		
	\end{definition}

    \begin{remark}
       For any $q\in (1,\infty]$, we have $L^{q}(\Omega)=(L^{p}(\Omega)'$ with  $\Omega$  in $\mathbb{R}^{3}$, i.e.
		\begin{equation*}
			L^{q}(\Omega)=\mathcal{M}(L^{p}(\Omega)\to L^{1}(\Omega))  \text{ with }\frac{1}{p}+\frac{1}{q}=1,
		\end{equation*}
	which means that there is some equivalence between the between the multiplier spaces and   Lebesgue spaces.
    \end{remark}
	
	\begin{remark}[Remark 1, \cite{CZ23}]\label{re2.2}
		For $\mathcal{M}^{3}(\dot{H}^{\frac{\theta}{\beta+\theta}}\to L^{2(\beta+\theta)})$, according to Sobolev embedding, we easily get $\dot{H}^{\frac{\theta}{\beta+\theta}}(\mathbb{R}^{3})\hookrightarrow L^{\frac{6(\beta+\theta)}{3\beta+\theta}}(\mathbb{R}^{3})$, since $\frac{1}{2(\beta+\theta)}=\frac{1}{\frac{6(\beta+\theta)}{3\beta+\theta}}+\frac{1}{\frac{6(\beta+\theta)}{3-3\beta-\theta}}$, we gain
		\begin{equation*}
			L^{\frac{6(\beta+\theta)}{3-3\beta-\theta}}(\mathbb{R}^{3})\hookrightarrow\mathcal{M}^{3}(\dot{H}^{\frac{\theta}{\beta+\theta}}(\mathbb{R}^{3})\to L^{2(\beta+\theta)}(\mathbb{R}^{3})).
		\end{equation*}
	\end{remark}

Using the interpolation method to obtain the boundedness of linear and nonlinear terms, we also need the following interpolation inequality:
	\begin{lemma}[Gagliardo-Nirenberg inequality,  ~Lemma 2.2, ~\cite{Wu24} or \cite{Nirenberg11}]\label{le2.3}
		Let $1\leq q,r<\infty$ and $m\leq k$.
		Suppose that $j$ and $\vartheta$ satisfy $m\leq j\leq k$, $0\leq \vartheta\leq 1$ and define $p\in [1,+\infty]$ by
		
		\begin{equation*}
			j- \frac{3}{p}=\left(m-\frac{3}{r}\right)\vartheta+\left(k-\frac{3}{q}\right)(1-\vartheta).
		\end{equation*}
		Then the inequality holds:
		\begin{equation*}
			\|\Lambda^{j} u\|_{L^{p}}\leq C	\|\Lambda^{m}u\|_{L^{r}}^{\vartheta}\|\Lambda^{k} u\|_{L^{q}}^{1-\vartheta}, ~ u\in W^{m,r}(\mathbb{R}^{3})\cap W^{k,q}(\mathbb{R}^{3}),
		\end{equation*}
		where constant $C\geq 0$. Here, when $p=\infty$, we require that $0<\vartheta<1$.
	\end{lemma}
	
	\begin{lemma}\label{le2.4}
		Let $k=0$, $j=\frac{\theta}{\beta+\theta}$, $m=\alpha$, $p=q=r=2$, and $\vartheta=\frac{\theta}{\alpha(\beta+\theta)}$ in Lemma \ref{le2.3}, and according to Remark \ref{re2.2},  we get the estimation
		\begin{eqnarray*}
			\||u|^{2}\|_{L^{2}} &=& \||u|^{\frac{2}{\beta+\theta}} \|^{\beta+\theta}_{L^{2(\beta+\theta)}}=\||u|^{\frac{2}{\beta+\theta}-2}u\cdot u\|^{\beta+\theta}_{L^{2(\beta+\theta)}} \\
			&\leq& \||u|^{\frac{2}{\beta+\theta}-2}u\|^{\beta+\theta}_{\mathcal{M}^{3}(\dot{H}^{\frac{\theta}{\beta+\theta}}\to L^{2(\beta+\theta)})}\|u\|^{\beta+\theta}_{\dot{H}^\frac{\theta}{\beta+\theta}} \\
			&\leq& \||u|^{\frac{2}{\beta+\theta}-2}u\|^{\beta+\theta}_{\mathcal{M}^{3}(\dot{H}^{\frac{\theta}{\beta+\theta}}\to L^{2(\beta+\theta)})}\|u\|^{\frac{\alpha(\beta+\theta)-\theta}{\alpha}}_{L^{2}}\|\Lambda^\alpha u\|^{\frac{\theta}{\alpha}}_{L^{2}}
		\end{eqnarray*}
		with
		\begin{equation*}
		      \begin{cases}
		      	\frac{1}{2(\beta+\theta)}=\frac{1}{\frac{6(\beta+\theta)}{3\beta+\theta}}+\frac{1}{\frac{6(\beta+\theta)}{3-3\beta-\theta}},\\
		      	\frac{\theta}{\beta+\theta}-\frac{3}{2}=(\alpha-\frac{3}{2})\frac{\theta}{\alpha(\beta+\theta)} + (-\frac{3}{2})\frac{\alpha(\beta+\theta)-\theta}{\alpha(\beta+\theta)}.
		      \end{cases}
		\end{equation*}
	\end{lemma}

Similar to how Yong and Jiu proved \cite[Lemma 3.1]{YJ09}, we obtain the fractional version of  \cite[Lemma 3.1]{YJ09} as follow:
	\begin{lemma}\label{le2.7}
		Suppose that $\begin{pmatrix}
			u(t)\\
			B(t)
		\end{pmatrix}$ is a weak solution to the fractional MHD on the interval $[0,T]$ with $0<T\le\infty$. Then $\begin{pmatrix}
			u(t)\\
			B(t)
		\end{pmatrix}$ is weakly continuous in the sense of $L^2$ after redefining the value of $\begin{pmatrix}
			u(t)\\
			B(t)
		\end{pmatrix}$  in a set of zero measure. Moreover, the following equalities hold:
		\begin{equation}\label{eq2.1}
			(u(t),\phi(t))=(u_0,\phi(0))-\int_{0}^{t}[(\Lambda^{\alpha} u,\Lambda^{\alpha}\phi)+((u\cdot\nabla)u,\phi)-((B\cdot\nabla)B,\phi) - (u,\phi_\tau)]{\rm d}\tau
		\end{equation}
		and
		\begin{equation}\label{eq2.2}
			(B(t),\phi(t))=(B_0,\phi(0))-\int_{0}^{t}[(\Lambda^{\alpha} B,\Lambda^{\alpha}\phi)+((u\cdot\nabla)B,\phi)-((B\cdot\nabla)u,\phi) - (B,\phi_\tau)]{\rm d}\tau,
		\end{equation}
where $u$ and  $B$ are in $L^{\infty}(0,T; L^{2}_{\sigma}(\mathbb{R}^{3}))\cap L^{2}(0,T; H^{\alpha}_{0}(\mathbb{R}^{3}))$ for any $t \in [0,T)$ and $\phi$ is in $C^{\infty}_c(0,T; C^{\infty}_{c,\sigma}(\mathbb{R}^{3}))$.
	\end{lemma}

	By applying the technique of symmetrization we transform the MHD equations \eqref{eq1.2} into the symmetric form, then we only need to consider the unknown variable $W$ that replaces $u$ and $B$.
	Let \begin{equation}\label{eq2.3}
		W(x,t)=\begin{pmatrix}
			u+B\\
			u-B
		\end{pmatrix}(x,t),
		~ \widetilde{\Lambda}^{\alpha} := \begin{pmatrix}
			\Lambda^{\alpha} & 0\\
			0 & \Lambda^{\alpha}
		\end{pmatrix}
		\text{ and }
		~ \widetilde{\nabla} := \begin{pmatrix}
			\nabla & 0\\
			0 & \nabla
		\end{pmatrix}.
	\end{equation}
		Then, we can easily get the following facts:
	\begin{description}
		\item[(\romannumeral1)]  $\|W\|^{2}_{L^{2}}=2\|u\|^{2}_{L^{2}}+2\|B\|^{2}_{L^{2}};$
		\item[(\romannumeral2)]  $\|W\|^{2}_{L^{2}(L^{2})}=2\|u\|^{2}_{L^{2}(L^{2})}+2\|B\|^{2}_{L^{2}(L^{2})};$
		\item[(\romannumeral3)]  $\|\widetilde{\Lambda}^{\alpha} W\|^{2}_{L^{2}(L^{2})}=2\|\Lambda^{\alpha} u\|^{2}_{L^{2}(L^{2})}+2\|\Lambda^{\alpha} B\|^{2}_{L^{2}(L^{2})};$
		\item[(\romannumeral4)]  $\widetilde{\nabla}\cdot W=0$.
	\end{description}
	Adding and subtracting between $(\ref{eq1.2})_1$ and $(\ref{eq1.2})_2$, then combining the nonlinear terms, we can obtain
	\begin{equation*}
		\begin{cases}
			(u+B)_{t} + (-\triangle)^{\alpha}(u+B) +[(u-B)\cdot \nabla)](u+B) + \nabla p=0,\quad &(x,t)\in \mathbb{R}^{3}\times (0,T),\\
			(u-B)_{t}+ (-\triangle)^{\alpha} (u-B) +[(u+B)\cdot \nabla)](u-B) + \nabla p=0,\quad &(x,t)\in \mathbb{R}^{3}\times (0,T).
		\end{cases}
	\end{equation*}
	We get
	\begin{equation*}
		W_{t}+ \begin{pmatrix}
			(-\triangle)^{\alpha} & 0\\
			0 & (-\triangle)^{\alpha}
		\end{pmatrix}
		W +\begin{pmatrix}
			0 & E\\
			E & 0
		\end{pmatrix} W \cdot
		\widetilde{\nabla} W +
		\widetilde{\nabla}
		\begin{pmatrix}
			p\\
			p
		\end{pmatrix} = 0,\quad (x,t)\in \mathbb{R}^{3}\times (0,T),
	\end{equation*}
	where $E$ is a third-order identity matrix.

	\begin{remark}\label{re2.9}
		Under the condition of Lemma \ref{le2.7}, we integrate \eqref{eq2.1} and \eqref{eq2.2} into a new equation by the technique of symmetrization:
		\begin{equation}\label{eq2.4}
			\begin{split}
				&(W, \Phi(t))=(W_0, \Phi(0))
				-\int_{0}^{t}\left[(\widetilde{\Lambda}^{\alpha} W,\widetilde{\Lambda}^{\alpha} \Phi)+\left( \begin{pmatrix}
					0 & E\\
					E & 0
				\end{pmatrix}W \cdot \widetilde{\nabla} W, \Phi\right )-(W, \Phi_\tau)\right]{\rm d}\tau
			\end{split}
		\end{equation}
		with
		\begin{equation*}
			W_0(x):= \begin{pmatrix}
				u_0 + B_0\\
				u_0 - B_0
			\end{pmatrix}(x)  \text{ and }
			\Phi(x,t) :=\begin{pmatrix}
				\phi(x,t)\\
				\phi(x,t)
			\end{pmatrix},
		\end{equation*}  $\phi \in C^{\infty}_c(0,T;C^{\infty}_{c,\sigma}(\mathbb{R}^{3})).$
	\end{remark}
	In fact, adding and subtracting between \eqref{eq2.1} and \eqref{eq2.2} separately, then combining the nonlinear terms, we have
	\begin{equation*}
		\begin{split}
			&(u+B,\phi(t))=(u_0 + B_0,\phi(0))\\
			&-\int_{0}^{t}\left[(\Lambda^{\alpha} (u+B),\Lambda^{\alpha}\phi)+\big( ((u-B)\cdot\nabla)(u+B),\phi\big )- (u+B,\phi_\tau)\right]{\rm d}\tau
		\end{split}
	\end{equation*}
	and
	\begin{equation*}
		\begin{split}
			&(u-B,\phi(t))=(u_0 - B_0,\phi(0))\\
			&-\int_{0}^{t}\left[(\Lambda^{\alpha} (u-B),\Lambda^{\alpha}\phi)+\big( ((u+B)\cdot\nabla)(u-B),\phi\big ) - (u-B,\phi_\tau) \right]{\rm d}\tau.
		\end{split}
	\end{equation*}
	This together with \eqref{eq2.3} implies \eqref{eq2.4}.

\section{The proof of Theorem 1.3
 }\label{Sec3}
	
	\begin{proof}[Proof of Theorem \ref{th1.1}]	
					
		Based on above symmetrization operation, in order to obtain the energy equality \eqref{eq1.3},  we only need to prove the equality containing $W$, that is,
		\begin{equation}\label{eq3.1}
			\|W(t_0)\|^{2}_{L^{2}}+2\int_{0}^{t_0}\|\widetilde{\Lambda}^{\alpha} W(\tau)\|^{2}_{L^{2}}{\rm d}\tau=\|W_0\|^{2}_{L^{2}}.
		\end{equation}

		The proof is divided into three steps.
		
		Step 1. Choosing suitable test function in \eqref{eq2.3}.\\
		Noting that  $C^{\infty}_{c,\sigma}(\mathbb{R}^{3})~$is dense in $H^{\alpha}_{0,\sigma}(\mathbb{R}^{3})$, we choose a sequence of $W^{j}\in C^{\infty}_c(0,T;C^{\infty}_{c,\sigma}(\mathbb{R}^{3}))$ converging to$\;W\;$in $L^{2}(0,T;L^{2}_{\sigma}(\mathbb{R}^{3}))\cap L^{2}(0,T;H^{\alpha}_{0,\sigma}(\mathbb{R}^{3}))$, meanwhile, we have
		\begin{equation*}
			W^{j}_\varepsilon(t)=J_{\varepsilon}[W^{j}](t)=\int_{0}^{t_0}\eta_\varepsilon(t-\sigma)W^{j}(\sigma){\rm d}\sigma
		\end{equation*}
		for any fixed $t_0<T$.
		
		Taking $W^{j}_\varepsilon\;$as test function in \eqref{eq2.3} and replacing $t\;$with fixed $t_0$, we get
		\begin{eqnarray}\label{eq3.2}
			& &\left(W(t_0),W^{j}_\varepsilon(t_0)\right)\nonumber\\
			&=&\int_{0}^{t_0}\eta_\varepsilon(t_0-\sigma)(W(t_0),W^{j}(\sigma)){\rm d}\sigma \nonumber\\
			&=&\int_{0}^{t_0}\eta_\varepsilon(-\sigma)(W_0,W^{j}(\sigma)){\rm d}\sigma \nonumber\\
			& & -\int_{0}^{t_0}\int_{0}^{t_0}\eta_\varepsilon(\tau-\sigma)\left[(\widetilde{\Lambda}^{\alpha} W(\tau),\widetilde{\Lambda}^{\alpha} W^{j}(\sigma))+\left( \begin{pmatrix}
				0 & E\\
				E & 0
			\end{pmatrix}W(\tau) \cdot \widetilde{\nabla} W(\tau),W^{j}(\sigma)\right )\right]{\rm d}\sigma {\rm d}\tau\nonumber\\
			& & +\int_{0}^{t_0}\int_{0}^{t_0}\frac{d}{d\tau}[\eta_\varepsilon(\tau-\sigma)](W(\tau),W^{j}(\sigma)){\rm d}\sigma {\rm d}\tau.
		\end{eqnarray}

		Step 2. Fixing $\varepsilon$ and sending $j\;$to infinity in \eqref{eq3.2}, then, we're going to estimate each of these terms in the above equation, starting with the first term on the left side of \eqref{eq3.2},
		\begin{eqnarray*}
			& & 	\left|\int_{0}^{t_0}\eta_\varepsilon(t_0-\sigma)(W(t_0),W^{j}(\sigma)-W(\sigma)){\rm d}\sigma\right| \\
			&\leq~~& \int_{0}^{t_0}\eta_\varepsilon(t_0-\sigma)\|W(t_0)\|_{L^{2}}\|W^{j}(\sigma)-W(\sigma)\|_{L^{2}}{\rm d}\sigma\quad(\text{ H\"{o}lder's inequality})\\
			&\lesssim_{\varepsilon,t_0}& \|W\|_{L^{\infty}(L^{2})}\|W^{j}-W\|_{L^{2}(L^{2})}\quad(\text{ H\"{o}lder's inequality})
		\end{eqnarray*}

		$\to0\quad as~j\to\infty,$
		we have
		\begin{equation*}
			\int_{0}^{t_0}\eta_\varepsilon(t_0-\sigma)(W(t_0),W^{j}(\sigma)){\rm d}\sigma\to\;\int_{0}^{t_0}\eta_\varepsilon(t_0-\sigma)(W(t_0),W(\sigma)){\rm d}\sigma~as~j\to \infty.
		\end{equation*}
		Similarly,
		\begin{equation*}
			\int_{0}^{t_0}\eta_\varepsilon(-\sigma)(W_0,W^{j}(\sigma)){\rm d}\sigma\to\;\int_{0}^{t_0}\eta_\varepsilon(-\sigma)(W_0,W(\sigma)){\rm d}\sigma~as~j\to \infty.
		\end{equation*}
		
		Since
		\begin{eqnarray*}
			& & \bigg|\int_{0}^{t_0}\int_{0}^{t_0}\eta_\varepsilon(\tau-\sigma)(\widetilde{\Lambda}^{\alpha} W(\tau),\widetilde{\Lambda}^{\alpha} W^{j}(\sigma)-\widetilde{\Lambda}^{\alpha} W(\sigma)){\rm d}\sigma {\rm d}\tau\bigg| \\
			&\lesssim_{\varepsilon}& \bigg|\int_{0}^{t_0}\int_{0}^{t_0}(\widetilde{\Lambda}^{\alpha} W(\tau),\widetilde{\Lambda}^{\alpha} W^{j}(\sigma)-\widetilde{\Lambda}^{\alpha} W(\sigma)){\rm d}\sigma{\rm d}\tau\bigg| \\
			&\lesssim_{\varepsilon}&  \int_{0}^{t_0}\int_{0}^{t_0}\|\widetilde{\Lambda}^{\alpha} W(\tau)\|_{L^{2}}\|\widetilde{\Lambda}^{\alpha} W^{j}(\sigma)-\widetilde{\Lambda}^{\alpha} W(\sigma)\|_{L^{2}}{\rm d}\sigma{\rm d}\tau \quad(\text{ H\"{o}lder's inequality})\\
			&=& C(\varepsilon)\int_{0}^{t_0}\|\widetilde{\Lambda}^{\alpha} W(\tau)\|_{L^{2}}{\rm d}\tau\int_{0}^{t_0}\|\widetilde{\Lambda}^{\alpha} W^{j}(\sigma)-\widetilde{\Lambda}^{\alpha} W(\sigma)\|_{L^{2}}{\rm d}\sigma \\
			&\lesssim_{\varepsilon}&\|\widetilde{\Lambda}^{\alpha} W\|_{L^{2}(L^{2})}\|\widetilde{\Lambda}^{\alpha} (W^{j}-W)\|_{L^{2}(L^{2})}
			\to\;0  \text{ ~ as~ }j\to\;\infty,\quad(\text{ H\"{o}lder's inequality})
		\end{eqnarray*}
		
		we get
		\begin{equation*}
			\int_{0}^{t_0}\int_{0}^{t_0}\eta_\varepsilon(\tau-\sigma)(\widetilde{\Lambda}^{\alpha} W(\tau),\widetilde{\Lambda}^{\alpha} W^{j}(\sigma)){\rm d}\sigma {\rm d}\tau\to\;\int_{0}^{t_0}\int_{0}^{t_0}\eta_\varepsilon(\tau-\sigma)(\widetilde{\Lambda}^{\alpha} W(\tau),\widetilde{\Lambda}^{\alpha} W(\sigma)){\rm d}\sigma {\rm d}\tau
		\end{equation*}
		as $j\to\;\infty$.
		
		Similarly, we have
		\begin{equation*}
			\int_{0}^{t_0}\int_{0}^{t_0}\frac{d}{d\tau}[\eta_\varepsilon(\tau-\sigma)](W(\tau),W^{j}(\sigma)){\rm d}\sigma {\rm d}\tau\to\;\int_{0}^{t_0}\int_{0}^{t_0}\frac{d}{d\tau}[\eta_\varepsilon(\tau-\sigma)](W(\tau),W(\sigma)){\rm d}\sigma {\rm d}\tau
		\end{equation*}
		as $j\to\;\infty$.
		
		Before estimating the nonlinear term $\begin{pmatrix}
			0&E\\
			E&0
		\end{pmatrix}W\cdot \widetilde{\nabla}W$, we have the following facts:
		\begin{equation*}
			\begin{split}
				&\big((u-B)\cdot\nabla(u+B),(u^{j}+B^{j})-(u+B)\big)=-\big((u-B)\cdot\nabla((u^{j}+B^{j})-(u+B)),u+B\big)\\
				&\text{and}\\
				&\big((u+B)\cdot\nabla(u-B),(u^{j}-B^{j})-(u-B)\big)=-\big((u+B)\cdot\nabla((u^{j}-B^{j})-(u-B)),u-B\big).
			\end{split}
		\end{equation*}
		Thus, we obtain
		\begin{equation*}
			\left( \begin{pmatrix}
				0 & E\\
				E & 0
			\end{pmatrix}W \cdot \widetilde{\nabla} W,W^{j}-W \right)= - \left( \begin{pmatrix}
				0 & E\\
				E & 0
			\end{pmatrix}W \cdot \widetilde{\nabla} (W^{j}-W),W \right).
		\end{equation*}
		Since
		\begin{equation*}
			W=\begin{pmatrix}
				u+B\\
				0
			\end{pmatrix}+\begin{pmatrix}
				0\\
				u-B
			\end{pmatrix}
		\end{equation*}
		and
		\begin{equation*}
			|W|\leq
			\begin{vmatrix}
				u+B\\
				0
			\end{vmatrix}+\begin{vmatrix}
				0\\
				u-B
			\end{vmatrix} \lesssim |u|+|B|,
		\end{equation*}
		we obtain
		\begin{equation*}
			|W|^{\frac{2}{\beta+\theta}-1} \lesssim |u|^{\frac{2}{\beta+\theta}-1}+	|B|^{\frac{2}{\beta+\theta}-1}.
		\end{equation*}
		Then we  get
		\begin{eqnarray*}
			& & \||W|^{\frac{2}{\beta+\theta}-2}W\|_{L^{\frac{2\alpha(\beta+\theta)}{2\alpha-\theta}}\mathcal{M}^{3}(\dot{H}^{\frac{\theta}{\beta+\theta}}\to L^{2(\beta+\theta)})}\\
			&\lesssim&  \||u|^{\frac{2}{\beta+\theta}-2}u\|_{L^{\frac{2\alpha(\beta+\theta)}{2\alpha-\theta}}\mathcal{M}^{3}(\dot{H}^{\frac{\theta}{\beta+\theta}}\to L^{2(\beta+\theta)})}+\||B|^{\frac{2}{\beta+\theta}-2}B\|_{L^{\frac{2\alpha(\beta+\theta)}{2\alpha-\theta}}\mathcal{M}^{3}(\dot{H}^{\frac{\theta}{\beta+\theta}}\to L^{2(\beta+\theta)})}.
		\end{eqnarray*}
		
		Now continue to estimate the nonlinear term $\begin{pmatrix}
			0&E\\
			E&0
		\end{pmatrix}W\cdot \widetilde{\nabla}W$, we have
		
		\begin{eqnarray*}
			& & \left|\int_{0}^{t_0}\int_{0}^{t_0}\eta_\varepsilon(\tau-\sigma)\left( \begin{pmatrix}
				0 & E\\
				E & 0
			\end{pmatrix}W(\tau) \cdot \widetilde{\nabla} W(\tau),W^{j}(\sigma)-W(\sigma)\right){\rm d}\sigma {\rm d}\tau\right|\\
			&\lesssim_{\varepsilon}& \left|\int_{0}^{t_0}\int_{0}^{t_0}\left( \begin{pmatrix}
				0 & E\\
				E & 0
			\end{pmatrix}W(\tau) \cdot \widetilde{\nabla} (W^{j}(\sigma)-W(\sigma)),W(\tau)\right){\rm d}\sigma {\rm d}\tau\right|\\
			&\lesssim_{\varepsilon}&	\int_{0}^{t_0}\int_{0}^{t_0}\|2|u|^{2}(\tau)-2|B|^{2}(\tau) \|_{L^{2}}\|\widetilde{\nabla} (W^{j}(\sigma)-W(\sigma))\|_{L^{2}}\;{\rm d}\sigma {\rm d}\tau \quad(\text{ H\"{o}lder's inequality})\\
			&\lesssim_{\varepsilon}&\int_{0}^{t_0}\int_{0}^{t_0}\|2|u|^{2}(\tau)+2|B|^{2}(\tau) \|_{L^{2}}\|\widetilde{\nabla} (W^{j}(\sigma)-W(\sigma))\|_{L^{2}}\;{\rm d}\sigma {\rm d}\tau \\
			&=& C(\varepsilon)\int_{0}^{t_0}\||W|^{2}(\tau)\|_{L^{2}}{\rm d}\tau\int_{0}^{t_0}\|\widetilde{\nabla} (W^{j}(\sigma)-W(\sigma))\|_{L^{2}}{\rm d}\sigma  \\
			&\lesssim_{\varepsilon}&\int_{0}^{t_0}\||W|^{\frac{2}{\beta+\theta}-2}W(\tau)\|^{\beta+\theta}_{\mathcal{M}^{3}(\dot{H}^{\frac{\theta}{\beta+\theta}}\to L^{2(\beta+\theta)})}\|W(\tau)\|^{\frac{\alpha(\beta+\theta)-\theta}{\alpha}}_{L^{2}}\|\widetilde{\Lambda}^\alpha W(\tau)\|^{\frac{\theta}{\alpha}}_{L^{2}}{\rm d}\tau\\
			& &\int_{0}^{t_0}\|\widetilde{\Lambda}^{\alpha} (W^{j}-W)(\sigma)\|^{\frac{1}{\alpha}}_{L^{2}}\|W^{j}(\sigma)-W(\sigma)\|^{\frac{\alpha-1}{\alpha}}_{L^{2}}{\rm d}\sigma\quad(\text{Gagliardo-Nirenberg inequality})\\
		    &\lesssim_{\varepsilon,t_0}& \left\|\||W|^{\frac{2}{\beta+\theta}-2}W\|^{\beta+\theta}_{\mathcal{M}^{3}(\dot{H}^{\frac{\theta}{\beta+\theta}}\to L^{2(\beta+\theta)})}\right\|_{L^{\frac{2\alpha}{2\alpha-\theta}}}\|W\|^{\frac{\alpha(\beta+\theta)-\theta}{\alpha}}_{L^{\infty}(L^{2})}\\
		    & & \left\|\|\widetilde{\Lambda}^{\alpha}W\|^{\frac{\theta}{\alpha}}_{L^{2}}\right\|_{L^{\frac{2\alpha}{\theta}}}\left\| \|\widetilde{\Lambda}^{\alpha} (W^{j}-W)\|^{\frac{1}{\alpha}}_{L^{2}}\right\|_{L^{2\alpha}} \left\|\|W^{j}-W\|^{\frac{\alpha-1}{\alpha}}_{L^{2}}\right\|_{L^{\frac{2\alpha}{\alpha-1}}}(\text{ H\"{o}lder's inequality}) \text{~\uppercase\expandafter{\romannumeral 1}}\\
			&=&C(\varepsilon,t_0)\||W|^{\frac{2}{\beta+\theta}-2}W\|^{\beta+\theta}_{L^{\frac{2\alpha(\beta+\theta)}{2\alpha-\theta}}\mathcal{M}^{3}(\dot{H}^{\frac{\beta}{\beta+\theta}}\to L^{2(\beta+\theta)})}\|W\|^{\frac{\alpha(\beta+\theta)-\theta}{\alpha}}_{L^{\infty}(L^{2})}\|\widetilde{\Lambda}^{\alpha}W\|^{\frac{\theta}{\alpha}}_{L^{2}(L^{2})}   \\
			& & \|\widetilde{\Lambda}^{\alpha} (W^{j}-W)\|^{\frac{1}{\alpha}}_{L^{2}(L^{2})} \|W^{j}-W\|^{\frac{\alpha-1}{\alpha}}_{L^{2}(L^{2})}\;\to0
		\end{eqnarray*}
		as$~j~\to \infty$, where the exponents at \uppercase\expandafter{\romannumeral 1} satisfy the following relation
		\begin{equation*}
			\begin{cases}
				1-\frac{3}{2}=(\alpha-\frac{3}{2})\frac{1}{\alpha}+(-\frac{3}{2})(1-\frac{1}{\alpha}),\\
				1=\frac{1}{\frac{2\alpha}{2\alpha-\theta}}+\frac{1}{\frac{2\alpha}{\theta}},\\
				\frac{1}{2}=\frac{1}{2\alpha}+\frac{1}{\frac{2\alpha}{\alpha-1}}.
			\end{cases}
		\end{equation*}
		Thus,
		\begin{eqnarray*}
			& & \int_{0}^{t_0}\int_{0}^{t_0}\eta_\varepsilon(\tau-\sigma)\left( \begin{pmatrix}
				0 & E\\
				E & 0
			\end{pmatrix}W(\tau) \cdot \widetilde{\nabla} W(\tau),W^{j}(\sigma)\right ){\rm d}\sigma {\rm d}\tau \\
			&\to& \int_{0}^{t_0}\int_{0}^{t_0}\eta_\varepsilon(\tau-\sigma)\left( \begin{pmatrix}
				0 & E\\
				E & 0
			\end{pmatrix}W(\tau) \cdot \widetilde{\nabla} W(\tau),W(\sigma)\right ){\rm d}\sigma {\rm d}\tau,\\
          \text{as~}j&\to& \infty.
		\end{eqnarray*}
		
		Finally, we obtain
		\begin{eqnarray}\label{eq3.3}
			& &\int_{0}^{t_0}\eta_\varepsilon(t_0-\sigma)(W(t_0),W(\sigma)){\rm d}\sigma\nonumber\\
			&=& \int_{0}^{t_0}\eta_\varepsilon(-\sigma)(W_0,W(\sigma))){\rm d}\sigma\nonumber\\
			& & -\int_{0}^{t_0}\int_{0}^{t_0}\eta_\varepsilon(\tau-\sigma)\left[(\widetilde{\Lambda}^{\alpha} W(\tau),\widetilde{\Lambda}^{\alpha} W(\sigma))\nonumber+\left( \begin{pmatrix}
				0 & E\\
				E & 0
			\end{pmatrix}W(\tau) \cdot \widetilde{\nabla} W(\tau),W(\sigma)\right )\right]{\rm d}\sigma {\rm d}\tau \nonumber\\
			& & +\int_{0}^{t_0}\int_{0}^{t_0}\frac{d}{d\tau}[\eta_\varepsilon(\tau-\sigma)](W(\tau),W(\sigma)){\rm d}\sigma {\rm d}\tau.
		\end{eqnarray}
		Step 3. Sending $\varepsilon$ to zero in \eqref{eq3.3}, and $\varepsilon<t_0$. For the term on the left of (\ref{eq3.3}), we have
		\begin{equation*}
			\int_{0}^{t_0}\eta_\varepsilon(t_0-\sigma)(W(t_0),W(\sigma)){\rm d}\sigma
			=\int_{0}^{\varepsilon}\eta_\varepsilon(\tau)(W(t_0),W(t_0-\tau))){\rm d}\tau.
		\end{equation*}
		Since $u$ and $B$ are weakly continuous in the sense of $L^{2}$, $W$ is weakly continuous in the sense of $L^{2}$, then
		\begin{equation*}
			(W(t_0),W(t_0-\tau))=\|W(t_0)\|^{2}_{L^{2}}+\widetilde{\beta}(\tau),
		\end{equation*}
		where $\widetilde{\beta}(\tau)\to\;0$ as $\sigma\to\;0$. Hence, as $\varepsilon\to\;0$ and mollifier is even, we have
		\begin{equation*}
			\int_{0}^{t_0}\eta_\varepsilon(t_0-\sigma)(W(t_0),W(\sigma)){\rm d}\sigma
			=\int_{0}^{\varepsilon}\eta_\varepsilon(\tau)(\|W(t_0)\|^{2}_{L^{2}}+\widetilde{\beta}(\tau)){\rm d}\tau
			\to\;\frac{1}{2}\|W(t_0)\|^{2}_{L^{2}}.
		\end{equation*}
		Similarly,
		\begin{equation*}
			\int_{0}^{t_0}\eta_\varepsilon(-\sigma)(W_0,W(\sigma)){\rm d}\sigma\to\;\frac{1}{2}\|W_0\|^{2}_{L^{2}}
		\end{equation*}as $\varepsilon\to 0$.\\
		For the second term on the right side of \eqref{eq3.3}, according to Lemma \ref{Lem2.0}, we obtain
		\begin{eqnarray*}
			& & \left|\int_{0}^{t_0}\int_{0}^{t_0}\eta_\varepsilon(\tau-\sigma)(\widetilde{\Lambda}^{\alpha} W(\tau), \widetilde{\Lambda}^{\alpha} W(\sigma)){\rm d}\sigma {\rm d}\tau-\int_{0}^{t_0}(\widetilde{\Lambda}^{\alpha} W(\tau), \widetilde{\Lambda}^{\alpha} W(\tau)){\rm d}\tau\right| \\
			&=& \bigg|\int_{0}^{t_0}(\widetilde{\Lambda}^{\alpha} W(\tau), \widetilde{\Lambda}^{\alpha} W_\varepsilon(\tau) - \widetilde{\Lambda}^{\alpha} W(\tau)){\rm d}\tau\bigg|\\
			&\le&\;\int_{0}^{t_0}\|\widetilde{\Lambda}^{\alpha} W(\tau)\|_{L^{2}}\|\widetilde{\Lambda}^{\alpha} W_\varepsilon(\tau) - \widetilde{\Lambda}^{\alpha} W(\tau)\|_{L^{2}}{\rm d}\tau \quad(\text{ H\"{o}lder's inequality})\\
			&\leq& \|\widetilde{\Lambda}^{\alpha} W\|_{L^{2}(L^{2})}\|\widetilde{\Lambda}^{\alpha} (W_\varepsilon-W)\|_{L^{2}(L^{2})}\to 0\quad(\text{ H\"{o}lder's inequality})
		\end{eqnarray*}
		as $\varepsilon\to 0$.
		
		Thus,
		\begin{equation*}
			\int_{0}^{t_0}\int_{0}^{t_0}\eta_\varepsilon(\tau-\sigma)(\widetilde{\Lambda}^{\alpha} W(\tau), \widetilde{\Lambda}^{\alpha} W(\sigma)){\rm d}\sigma {\rm d}\tau\to\;\int_{0}^{t_0}(\widetilde{\Lambda}^{\alpha}  W(\tau), \widetilde{\Lambda}^{\alpha} W(\tau)){\rm d}\tau\
		\end{equation*} as $\varepsilon\to 0$.\\
		
		For the nonlinear term $\begin{pmatrix}
			0&E\\
			E&0
		\end{pmatrix}W\cdot \widetilde{\nabla}W$, we have
		\begin{eqnarray*}
			& & \bigg|\int_{0}^{t_0}\int_{0}^{t_0}\eta_\varepsilon(\tau-\sigma)\left( \begin{pmatrix}
				0 & E\\
				E & 0
			\end{pmatrix}W(\tau) \cdot \widetilde{\nabla} W(\tau),W(\sigma)\right ){\rm d}\sigma {\rm d}\tau\\
			& & -\int_{0}^{t_0}\left( \begin{pmatrix}
				0 & E\\
				E & 0
			\end{pmatrix}W(\tau) \cdot \widetilde{\nabla} W(\tau),W(\tau)\right ) {\rm d}\tau\bigg|\\
			&=&	\left|\int_{0}^{t_0}\left( \begin{pmatrix}
				0 & E\\
				E & 0
			\end{pmatrix}W(\tau) \cdot \widetilde{\nabla}W(\tau), (W_\varepsilon(\tau)-W(\tau)\right){\rm d}\tau \right|\\
			&=& \left|\int_{0}^{t_0}\left( \begin{pmatrix}
				0 & E\\
				E & 0
			\end{pmatrix}W(\tau) \cdot \widetilde{\nabla} (W_\varepsilon(\tau)-W(\tau)),W(\tau)\right ) {\rm d}\tau\right|\\
			&\leq& 	\int_{0}^{t_0}\|2|u|^{2}(\tau)-2|B|^{2}(\tau) \|_{L^{2}}\|\widetilde{\nabla} (W_\varepsilon(\tau)-W(\tau))\|_{L^{2}}\; {\rm d}\tau\quad(\text{ H\"{o}lder's inequality}) \\
			&\leq&\int_{0}^{t_0}\|2|u|^{2}(\tau)+2|B|^{2}(\tau) \|_{L^{2}}\|\widetilde{\nabla} (W_\varepsilon(\tau)-W(\tau))\|_{L^{2}}\; {\rm d}\tau\\
			&=& \int_{0}^{t_0}\||W|^{2}(\tau) \|_{L^{2}}\|\widetilde{\nabla} (W_\varepsilon(\tau)-W(\tau))\|_{L^{2}}\; {\rm d}\tau\\
			&\leq&
			\int_{0}^{t_0}  \||W|^{\frac{2}{\beta+\theta}-2}W(\tau)\|^{\beta+\theta}_{\mathcal{M}^{3}(\dot{H}^{\frac{\theta}{\beta+\theta}}\to L^{2(\beta+\theta)})}\|W(\tau)\|^{\frac{\alpha(\beta+\theta)-\theta}{\alpha}}_{L^{2}}\|\widetilde{\Lambda}^{\alpha}W(\tau)\|^{\frac{\theta}{\alpha}}_{L^{2}}\\
			& &\|\widetilde{\Lambda}^{\alpha} (W^{j}_\varepsilon-W)(\tau)\|^{\frac{1}{\alpha}}_{L^{2}}\|W_\varepsilon(\tau)-W(\tau)\|^{\frac{\alpha-1}{\alpha}}_{L^{2}}{\rm d}\tau \quad(\text{Gagliardo-Nirenberg inequality})\\
			&\leq& \left\|\||W|^{\frac{2}{\beta+\theta}-2}W\|^{\beta+\theta}_{\mathcal{M}^{3}(\dot{H}^{\frac{\theta}{\beta+\theta}}\to L^{2(\beta+\theta)})}\right\|_{L^{\frac{2\alpha}{\alpha-\theta}}}\|W\|^{\frac{\alpha(\beta+\theta)-\theta}{\alpha}}_{L^{\infty}(L^{2})}\left\|\|\widetilde{\Lambda}^{\alpha}W\|^{\frac{\theta}{\alpha}}_{L^{2}}\right\|_{L^{\frac{2\alpha}{\theta}}}\\
			& &\left\| \|\widetilde{\Lambda}^{\alpha} (W_\varepsilon-W)\|^{\frac{1}{\alpha}}_{L^{2}}\right\|_{L^{2\alpha}} \left\|\|W_\varepsilon-W\|^{\frac{\alpha-1}{\alpha}}_{L^{2}}\right\|_{L^{\frac{2\alpha}{\alpha-1}}}    \quad(\text{ H\"{o}lder's inequality})\text{~\uppercase\expandafter{\romannumeral 2}}\\
			&=&\||W|^{\frac{2}{\beta+\theta}-2}W\|^{\beta+\theta}_{L^{\frac{2\alpha(\beta+\theta)}{\alpha-\theta}}\mathcal{M}^{3}(\dot{H}^{\frac{\beta}{\beta+\theta}}\to L^{2(\beta+\theta)})}\|W\|^{\frac{\alpha(\beta+\theta)-\theta}{\alpha}}_{L^{\infty}(L^{2})}\|\widetilde{\Lambda}^{\alpha}W\|^{\frac{\theta}{\alpha}}_{L^{2}(\dot{H}^{1})}   \\
			& & \|\widetilde{\Lambda}^{\alpha} (W_\varepsilon-W)\|^{\frac{1}{\alpha}}_{L^{2}(L^{2})} \|W_\varepsilon-W\|^{\frac{\alpha-1}{\alpha}}_{L^{2}(L^{2})}~\to0
		\end{eqnarray*}
		as $~\varepsilon~\to 0$, $j\to \infty$,  where the exponents at \uppercase\expandafter{\romannumeral 2} satisfy the following relation
		\begin{equation*}
			\begin{cases}
				1-\frac{3}{2}=(\alpha-\frac{3}{2})\frac{1}{\alpha}+(-\frac{3}{2})(1-\frac{1}{\alpha}),\\
				1=\frac{1}{\frac{2\alpha}{\alpha-\theta}}+\frac{1}{\frac{2\alpha}{\theta}}+\frac{1}{2\alpha}+\frac{1}{\frac{2\alpha}{\alpha-1}},\\
				\frac{1}{2}=\frac{1}{2\alpha}+\frac{1}{\frac{2\alpha}{\alpha-1}}.
			\end{cases}
		\end{equation*}
		Therefore, we have
		\begin{eqnarray*}
			& & \int_{0}^{t_0}\int_{0}^{t_0}\eta_\varepsilon(\tau-\sigma)\left( \begin{pmatrix}
				0 & E\\
				E & 0
			\end{pmatrix}W(\tau) \cdot \widetilde{\nabla} W(\tau),W(\sigma)\right ){\rm d}\sigma {\rm d}\tau \\
			&\to& \int_{0}^{t_0}\left( \begin{pmatrix}
				0 & E\\
				E & 0
			\end{pmatrix}W(\tau) \cdot \widetilde{\nabla} W(\tau),W(\tau)\right ) {\rm d}\tau
		\end{eqnarray*}
		as $~\varepsilon~\to 0$.
		
		Noting that
		\begin{equation*}
			\begin{split}
				&\bigg([(u-B)\cdot\nabla](u+B),u+B\bigg)=-\bigg([(u-B)\cdot\nabla](u+B),u+B\bigg)=0\\
\text{and}\\
				&\bigg([(u+B)\cdot\nabla](u-B),u-B\bigg)=-\bigg([(u+B)\cdot\nabla](u-B),u-B\bigg)=0,
			\end{split}
		\end{equation*}
		thus
		\begin{equation*}
			\left( \begin{pmatrix}
				0 & E\\
				E & 0
			\end{pmatrix}W(\tau) \cdot \widetilde{\nabla} W(\tau),W(\tau)\right )=0,
		\end{equation*}
		we have the fact that
		
		\begin{equation*}
			\int_{0}^{t_0}\int_{0}^{t_0}\eta_\varepsilon(\tau-\sigma)\left( \begin{pmatrix}
				0 & E\\
				E & 0
			\end{pmatrix}W(\tau) \cdot \widetilde{\nabla} W(\tau),W(\sigma)\right ){\rm d}\sigma {\rm d}\tau\to\;0\quad as~\varepsilon~\to 0.
		\end{equation*}
		For the last term of Equation \eqref{eq3.3},  given that $\eta_\varepsilon$ is an even function, we have
		\begin{eqnarray*}
			& & \int_{0}^{t_0}\int_{0}^{t_0}\frac{d}{d\tau}[\eta_\varepsilon(\tau-\sigma)](W(\tau),W(\sigma)){\rm d}\sigma {\rm d}\tau \\
			&=& -\int_{0}^{t_0}\int_{0}^{t_0}\frac{d}{d\sigma}[\eta_\varepsilon(\tau-\sigma)](W(\tau),W(\sigma)){\rm d}\sigma {\rm d}\tau \\
			&=& -\int_{0}^{t_0}\int_{0}^{t_0}\frac{d}{d\sigma}[\eta_\varepsilon(\sigma-\tau)](W(\tau),W(\sigma)){\rm d}\sigma {\rm d}\tau \\
			&=& -\int_{0}^{t_0}\int_{0}^{t_0}\frac{d}{d\sigma}[\eta_\varepsilon(\sigma-\tau)](W(\sigma),W(\tau)){\rm d}\tau {\rm d}\sigma \\
			&=& -\int_{0}^{t_0}\int_{0}^{t_0}\frac{d}{d\tau}[\eta_\varepsilon(\tau-\sigma)](W(\tau),W(\sigma)){\rm d}\sigma {\rm d}\tau.
		\end{eqnarray*}
		Therefore,
		\begin{equation*}
			\int_{0}^{t_0}\int_{0}^{t_0}\frac{d}{d\tau}[\eta_\varepsilon(\tau-\sigma)](W(\tau),W(\sigma)){\rm d}\sigma {\rm d}\tau=0.
		\end{equation*}
		Concluding the results and sending $\varepsilon\to\;0\;$ in (3.4), we get
		\begin{equation*}
			\frac{1}{2}\|W(t_0)\|^{2}_{L^{2}}=\frac{1}{2}\|W_0\|^{2}_{L^{2}}-\int_{0}^{t_0}\|\widetilde{\Lambda}^{\alpha} W(\tau)\|^{2}_{L^{2}}{\rm d}\tau,
		\end{equation*}
		i.e. \eqref{eq3.1}.
	\end{proof}
		\begin{remark}
		Taking $B=0$ in \eqref{eq2.4}, by using the previous method, we can get similar estimations about $u$, which provides the validity of Corollary \ref{co1.2}.
	\end{remark}
	\section{The proof of Theorem 1.6
and Theorem 1.9
}\label{Sec4}
	\begin{proof}[The sketch of proof of Theorem \ref{th1.3}]	
Since the estimations of  the linear terms in \eqref{eq3.2} are the same as Theorem \ref{th1.1},  we omit them. And we only give the estimation about the nonlinear term,
		For Theorem \ref{th1.3}, we have
\begin{eqnarray*}
	 & &\left|\int_{0}^{t_0}\left( \begin{pmatrix}
	 	0 & E\\
	 	E & 0
	 \end{pmatrix}W(\tau) \cdot \widetilde{\nabla} (W^{j}_\varepsilon(\tau)-W(\tau)),W(\tau)\right ) {\rm d}\tau\right|\\
	 &\lesssim&\int_{0}^{t_0}\left|\left( W(\tau) \cdot \widetilde{\nabla} (W^{j}_\varepsilon(\tau)-W(\tau)),W(\tau)\right )\right| {\rm d}\tau\\
	 &\leq&\int_{0}^{t_0}\left\||W|^{\gamma}~\widetilde{\nabla} (W^{j}_\varepsilon-W) (\tau)\right\|_{L^{\frac{1}{\gamma}}}\||W|^{2-\gamma}(\tau)\|_{L^{\frac{1}{1-\gamma}}}{\rm d}\tau\quad(\text{ H\"{o}lder's inequality})\\
	 &=&\int_{0}^{t_0}\left\|W~\left|\widetilde{\nabla} (W^{j}_\varepsilon-W)\right|^{\frac{1}{\gamma}} (\tau)\right\|^{\gamma}_{L^{1}} \|W(\tau)\|^{2-\gamma}_{L^{\frac{2-\gamma}{1-\gamma}}}{\rm d}\tau\\
	 &\leq&\int_{0}^{t_0} \left\|\left|\widetilde{\nabla} (W^{j}_\varepsilon-W)\right|^{\frac{1}{\gamma}} (\tau) \right\|^{\gamma}_{\mathcal{M}^{3}(\dot{H}^{\frac{1}{\beta+1}}\to L^{1})}\|W(\tau)\|^{\gamma}_{\dot{H}^{\frac{1}{\beta+1}}}\|W(\tau)\|^{2-\gamma}_{L^{\frac{2-\gamma}{1-\gamma}}}{\rm d}\tau\\
	 &\leq&\int_{0}^{t_0} \left\|\left|\widetilde{\nabla} (W^{j}_\varepsilon-W)\right|^{\frac{1}{\gamma}} (\tau) \right\|^{\gamma}_{\mathcal{M}^{3}(\dot{H}^{\frac{1}{\beta+1}}\to L^{1})}\|W(\tau)\|^{\gamma}_{\dot{H}^{\frac{1}{\beta+1}}}\|W(\tau)\|^{2-\gamma}_{L^{\frac{2-\gamma}{1-\gamma}}}{\rm d}\tau\\
	 &\lesssim&\int_{0}^{t_0} \left\|\left|\widetilde{\nabla} (W^{j}_\varepsilon-W)\right|^{\frac{1}{\gamma}} (\tau) \right\|^{\gamma}_{\mathcal{M}^{3}(\dot{H}^{\frac{1}{\beta+1}}\to L^{1})}\|\Lambda^{\alpha} W(\tau)\|^{\frac{\gamma}{\alpha(\beta+1)}}_{L^{2}}\|W(\tau)\|^{\frac{\alpha\gamma(\beta+1)-\gamma}{\alpha(\beta+1)}}_{L^{2}}\|W(\tau)\|^{2-\gamma}_{L^{\frac{2-\gamma}{1-\gamma}}}{\rm d}\tau\\
	 & &\quad(\text{Gagliardo-Nirenberg inequality})\\
	 &\leq&\left\| \left\|\left|\widetilde{\nabla} (W^{j}_\varepsilon-W)\right|^{\frac{1}{\gamma}} (\tau) \right\|^{\gamma}_{\mathcal{M}^{3}(\dot{H}^{\frac{1}{\beta+1}}\to L^{1})}  \right\|     _{L^{\frac{2\alpha(\beta+1)}{4\alpha\beta+4\alpha-4\alpha\beta\gamma-4\alpha\gamma-\gamma}}}\left\|\|\Lambda^{\alpha} W\|^{\frac{\gamma}{\alpha(\beta+1)}}_{L^{2}}\right\|_{L^{\frac{2\alpha(\beta+1)}{\gamma}}}\\
	 & &\|W\|^{\frac{\alpha\gamma(\beta+1)-\gamma}{\alpha(\beta+1)}}_{L^{\infty}L^{2}}~\left\|\|W\|^{2-\gamma}_{L^{\frac{2-\gamma}{1-\gamma}}}\right\|_{L^{\frac{1}{2\gamma-1}}}\quad(\text{ H\"{o}lder's inequality})\\
	 &=&\left\|\left|\widetilde{\nabla} (W^{j}_\varepsilon-W)\right|^{\frac{1}{\gamma}} \right\|^{\gamma}_{L^{\frac{2\alpha\gamma(\beta+1)}{4\alpha\beta+4\alpha-4\alpha\beta\gamma-4\alpha\gamma-\gamma}}\mathcal{M}^{3}(\dot{H}^{\frac{1}{\beta+1}}\to L^{1})}\|\Lambda^{\alpha} W\|^{\frac{\gamma}{\alpha(\gamma+1)}}_{L^{2}(L^{2})}\\
     & &\|W\|^{\frac{\alpha\gamma(\beta+1)-\gamma}{\alpha(\beta+1)}}_{L^{\infty}(L^{2})}~\|W\|^{2-\gamma}_{L^{\frac{2-\gamma}{2\gamma-1}}(L^{\frac{2-\gamma}{1-\gamma}})},
\end{eqnarray*}where  the exponents satisfy the following relation,
	\begin{equation*}
	\begin{cases}
		\frac{1}{\beta+1}-\frac{3}{2}=(\alpha-\frac{3}{2})\frac{1}{\alpha(\beta+1)} + (-\frac{3}{2})(1-\frac{1}{\alpha(\beta+1)}),\\
		1=\frac{1}{\frac{2\alpha(\beta+1)}{\gamma}}+\frac{1}{\frac{2\alpha(\beta+1)}{4\alpha\beta+4\alpha-4\alpha\beta\gamma-4\alpha\gamma-\gamma}}+\frac{1}{\frac{1}{2\gamma-1}}.
	\end{cases}
\end{equation*}
where $\gamma=\frac{7-3\theta}{10-5\theta}$. For $\|W\|_{L^{\frac{2-\gamma}{2\gamma-1}},L^{\frac{2-\gamma}{1-\gamma}}}$, taking $s=\frac{15\gamma-10}{5\gamma-3}$, we have
\begin{eqnarray*}
	& &\|W\|^{\frac{2\gamma-1}{2-\gamma}}_{L^{\frac{2-\gamma}{2\gamma-1}},L^{\frac{2-\gamma}{1-\gamma}}}\\
	&=&\int_{0}^{t_0}\|W(\tau)\|^{{\frac{2-\gamma}{2\gamma-1}}}_{L^{\frac{2-\gamma}{1-\gamma}}}{\rm d}\tau\\
	&\leq&\int_{0}^{t_0}\left\| |W|^{\frac{2-\gamma}{s-\gamma}-2}W\cdot W(\tau)\right\|^{\frac{s-\gamma}{2\gamma-1}}_{L^{\frac{s-\gamma}{1-\gamma}}}{\rm d}\tau\\
	&\leq&\int_{0}^{t_0}\left\| |W|^{\frac{2-\gamma}{s-\gamma}-2}W(\tau)\right\|^{\frac{s-\gamma}{2\gamma-1}}_{\mathcal{M}^{3}(\dot{H}^{1}\to L^{\frac{s-\gamma}{1-\gamma}})}\|W(\tau)\|^{\frac{s-\gamma}{2\gamma-1}}_{\dot{H}^{1}}{\rm d}\tau\\
	&\leq&\int_{0}^{t_0}\left\| |W|^{\frac{2-\gamma}{s-\gamma}-2}W(\tau)\right\|^{\frac{s-\gamma}{2\gamma-1}}_{\mathcal{M}^{3}(\dot{H}^{1}\to L^{\frac{s-\gamma}{1-\gamma}})}\|\Lambda^{\alpha}W(\tau)\|^{\frac{s-\gamma}{\alpha(2\gamma-1)}}_{L^{2}}\|W(\tau)\|^{\frac{(\alpha-1)(s-\gamma)}{\alpha(2\gamma-1)}}_{L^{2}}{\rm d}\tau\\
	& &\quad(\text{Gagliardo-Nirenberg inequality})\\
	&\leq&\left\|\left\| |W|^{\frac{2-\gamma}{s-\gamma}-2}W\right\|^{\frac{s-\gamma}{2\gamma-1}}_{\mathcal{M}^{3}(\dot{H}^{1}\to L^{\frac{s-\gamma}{1-\gamma}})}\right\|_{L^{\frac{2\alpha(2\gamma-1)}{2\alpha(2\gamma-1)-s+\gamma}}}\left\|\|\Lambda^{\alpha}W\|^{\frac{(s-\gamma)}{\alpha(2\gamma-1)}}_{L^{2}} \right\|_{L^{\frac{2\alpha(2\gamma-1)}{s-\gamma}}}\|W\|^{\frac{(\alpha-1)(s-\gamma)}{\alpha(2\gamma-1)}}_{L^{\infty}(L^{2})}\\
	& &\quad(\text{ H\"{o}lder's inequality})\\
	&=&\left\| |W|^{\frac{2-\gamma}{s-\gamma}-2}W\right\|^{\frac{s-\gamma}{2\gamma-1}}_{{L^{\frac{2\alpha(s-\gamma)}{2\alpha(2\gamma-1)-s+\gamma}}}\mathcal{M}^{3}(\dot{H}^{1}\to L^{\frac{s-\gamma}{1-\gamma}})}\|\Lambda^{\alpha}W\|^{\frac{(s-\gamma)}{\alpha(2\gamma-1)}}_{L^{2}(L^{2})}\|W\|^{\frac{(\alpha-1)(s-\gamma)}{\alpha(2\gamma-1)}}_{L^{\infty}(L^{2})},
\end{eqnarray*}where the exponents satisfy the following relation,
	\begin{equation*}
	\begin{cases}
		1-\frac{3}{2}=(\alpha-\frac{3}{2})\frac{1}{\alpha} + (-\frac{3}{2})(1-\frac{1}{\alpha}),\\
		1=\frac{1}{\frac{2\alpha(2\gamma-1)}{2\alpha(2\gamma-1)-s+\gamma}}+\frac{1}{\frac{2\alpha(2\gamma-1)}{s-\gamma}}.
	\end{cases}
\end{equation*}
\end{proof}

\begin{proof}[The sketch of proof of Theorem \ref{th1.5}]
For simplicity, we only need to estimate the nonlinear part in \eqref{eq3.2}.
 For Theorem \ref{th1.5}, we have
	\begin{eqnarray*}
		& &\left|\int_{0}^{t_0}\left( \begin{pmatrix}
	 	0 & E\\
	 	E & 0
	 \end{pmatrix}W(\tau) \cdot \widetilde{\nabla} W(\tau),W^{j}_\varepsilon(\tau)-W(\tau)\right ) {\rm d}\tau\right|\\
		&\leq&\int_{0}^{t_0}\left\|W^{j}_\varepsilon(\tau)-W (\tau)\right\|_{L^{\frac{12-6\theta}{3-\theta}}}\|W\cdot \widetilde{\nabla} W(\tau)\|_{L^{\frac{12-6\theta}{9-5\theta}}}{\rm d}\tau\quad(\text{ H\"{o}lder's inequality})\\
		&\leq&\int_{0}^{t_0}\left\|W^{j}_\varepsilon(\tau)-W (\tau)\right\|_{L^{\frac{12-6\theta}{3-\theta}}}\|\widetilde{\nabla} W\|_{\mathcal{M}^{3}(\dot{W}^{\frac{1+\theta}{1+\beta},\frac{6-3\theta}{3-\theta^{2}}}\to L^{\frac{12-6\theta}{9-5\theta}})}\|W\|_{\dot{W}^{\frac{1+\theta}{1+\beta},\frac{6-3\theta}{3-\theta^{2}}}}{\rm d}\tau\\
		&\lesssim& \int_{0}^{t_0}\left\|W^{j}_\varepsilon(\tau)-W (\tau)\right\|_{L^{\frac{12-6\theta}{3-\theta}}}\|\widetilde{\nabla} W\|_{\mathcal{M}^{3}(\dot{W}^{\frac{1+\theta}{1+\beta},\frac{6-3\theta}{3-\theta^{2}}}\to L^{\frac{12-6\theta}{9-5\theta}})}\|\widetilde{\Lambda}^{\alpha}W\|^{\frac{4-\theta-3\beta\theta+2\theta^{2}\beta}{2\alpha(1+\beta)(2-\theta)}}_{L^{2}}\\
		& &\|W\|^{1-\frac{4-\theta-3\beta\theta+2\theta^{2}\beta}{2\alpha(1+\beta)(2-\theta)}}_{L^{2}}{\rm d}\tau\quad(\text{Gagliardo-Nirenberg inequality})\\
		&=&C\int_{0}^{t_0}\left\|W^{j}_\varepsilon(\tau)-W (\tau)\right\|_{L^{\frac{12-6\theta}{3-\theta}}}\|\widetilde{\nabla} W\|_{\mathcal{M}^{3}(\dot{W}^{\frac{1+\theta}{1+\beta},\frac{6-3\theta}{3-\theta^{2}}}\to L^{\frac{12-6\theta}{9-5\theta}})}\|\widetilde{\Lambda}^{\alpha}W\|^{\frac{4-\theta-3\beta\theta+2\theta^{2}\beta}{\alpha(4-2\theta+4\beta-2\theta\beta)}}_{L^{2}}\\
		& &\|W\|^{1-\frac{4-\theta-3\beta\theta+2\theta^{2}\beta}{2\alpha(1+\beta)(2-\theta)}}_{L^{2}}{\rm d}\tau\quad ~~(\text{\uppercase\expandafter{\romannumeral 3}})\\
		&\lesssim& \int_{0}^{t_0}\left\|W^{j}_\varepsilon(\tau)-W (\tau)\right\|_{L^{\frac{12-6\theta}{3-\theta}}}\|\widetilde{\nabla} W\|_{\mathcal{M}^{3}(\dot{W}^{\frac{1+\theta}{1+\beta},\frac{6-3\theta}{3-\theta^{2}}}\to L^{\frac{12-6\theta}{9-5\theta}})}(\|\widetilde{\Lambda}^{\alpha}W\|_{L^{2}}+1)\\
		& &\|W\|^{1-\frac{4-\theta-3\beta\theta+2\theta^{2}\beta}{2\alpha(1+\beta)(2-\theta)}}_{L^{2}}{\rm d}\tau\\
		&\lesssim&\left\|W^{j}_\varepsilon-W \right\|_{L^{q}(L^{\frac{12-6\theta}{3-\theta}})}\|\widetilde{\nabla} W\|_{L^{p}\mathcal{M}^{3}(\dot{W}^{\frac{1+\theta}{1+\beta},\frac{6-3\theta}{3-\theta^{2}}}\to L^{\frac{12-6\theta}{9-5\theta}})}\left\|\|\widetilde{\Lambda}^{\alpha}W\|_{L^{2}}+1\right\|_{L^{2}}\\
		& &\|W\|^{1-\frac{4-\theta-3\beta\theta+2\theta^{2}\beta}{2\alpha(1+\beta)(2-\theta)}}_{L^{\infty}(L^{2})}\quad(\text{ H\"{o}lder's inequality})\\
		&\lesssim_{t_0}&\left\|W^{j}_\varepsilon-W \right\|_{L^{q}(L^{\frac{12-6\theta}{3-\theta}})}\|\widetilde{\nabla} W\|_{L^{p}\mathcal{M}^{3}(\dot{W}^{\frac{1+\theta}{1+\beta},\frac{6-3\theta}{3-\theta^{2}}}\to L^{\frac{12-6\theta}{9-5\theta}})}(\|\widetilde{\Lambda}^{\alpha}W\|_{L^{2}(L^{2})}+1)\\
		& &\|W\|^{1-\frac{4-\theta-3\beta\theta+2\theta^{2}\beta}{2\alpha(1+\beta)(2-\theta)}}_{L^{\infty}(L^{2})},
	\end{eqnarray*} where the exponents satisfy the following relation,
		\begin{equation*}
		\begin{cases}
			\frac{1+\theta}{1+\beta}-\frac{3}{\frac{6-3\theta}{3-\theta^{2}}} =  (\alpha-\frac{3}{2})\frac{4-\theta-3\beta\theta+2\theta^{2}\beta}{2\alpha(1+\beta)(2-\theta)} +  (-\frac{3}{2})(1-\frac{4-\theta-3\beta\theta+2\theta^{2}\beta}{2\alpha(1+\beta)(2-\theta)}),\\
			1=\frac{1}{2}+\frac{1}{q}+\frac{1}{p},~2\leq p<\infty.
		\end{cases}
	\end{equation*}
	
\end{proof}

\begin{remark}
  Taking $B=0$ in \eqref{eq2.4}, we can get similar estimations about $u$, which proves the validity of Corollary \ref{co1.4} and  Corollary \ref{co1.6}.
\end{remark}

Actually, we can replace the exponent of H\"{o}lder's inequality $	1=\frac{1}{2}+\frac{1}{q}+\frac{1}{p}$ with $	1=\frac{1}{\frac{2\alpha(4-2\theta+4\beta-2\theta\beta)}{4-\theta-3\beta\theta+2\theta^{2}\beta}}+\frac{1}{q'}+\frac{1}{p'}$ at (\text{\uppercase\expandafter{\romannumeral 3}}), where $\frac{2\alpha(4-2\theta+4\beta-2\theta\beta)}{4-\theta-3\beta\theta+2\theta^{2}\beta}>2\alpha\geq2$, $p'$ has a wider range of values than $p$, and Theorem \ref{th1.5} and Corollary \ref{co1.6} still valid.

\section*{Data availability statement}
No new data were created or analysed in this study.
\section*{Author contributions}
Yi Feng: Formal analysis (equal); Writing - original draft (equal).
Weihua Wang: Formal analysis (equal); Writing - original draft (equal); Writing - review \& editing (equal).

\section*{Declarations}
{\bf Conflict of interest}~~~~ The authors declare that they have no Conflict of interest.

%


	\end{document}